\documentclass[reqno,11pt]{amsart}
\usepackage{graphicx}
\usepackage{amsmath,amssymb}
\newtheorem{theorem}{Theorem}[section]
\newtheorem{proposition}[theorem]{Proposition}
\theoremstyle{definition}
\newtheorem{definition}[theorem]{Definition}

\numberwithin{equation}{section}
\begin{document}

\title[Dynamics of spatial logistic model]{Dynamics of spatial logistic model: finite systems}%
\author{Yuri Kozitsky}%
\address{Institute of Mathematics, Maria Curie-Sk{\l}odowska University, Lublin 20-031 Poland}%
\email{jkozi@hektor.umcs.lublin.pl}%

\thanks{Supported by the DFG through SFB 701: ``Spektrale
Strukturen und Topologische Methoden in der Mathematik", by the ZiF
Research Group ``Stochastic Dynamics: Mathematical Theory and
Applications" (Universit\"at Bielefeld), and by the European Commission under the project
STREVCOMS PIRSES-2013-612669.}%
\subjclass{82C22, 92D25, 60J80}%
\keywords{Individual-based model;  birth-and-death process; random point field}%

\begin{abstract}
The spatial logistic model is a system of point entities
(particles) in $\mathbb{R}^d$ which reproduce themselves at distant points (dispersal)
and die, also due to competition. The states of such systems are probability measures on the space of all locally finite particle configurations.
In this paper, we obtain the evolution of states  of `finite systems', that is, in the case where the initial state is supported on the subset of the configuration space consisting of finite configurations. The evolution is obtained as the global solution of
the corresponding Fokker-Planck equation  in the space of measures supported on the set of finite configurations. We also prove that  this evolution preserves the existence of exponential moments and the absolute continuity with respect to the Lebesgue-Poisson measure.
\end{abstract}
\maketitle

\section{Introduction}
\label{sec:1}

In a number of applications, one deals with large systems of
interacting entities evolving in time and distributed over an unbounded
continuous habitat. In simple models of such systems, each entity is
completely characterized by its spatial location  $x\in
\mathbb{R}^d$, $d\geq 1$. The entities form discrete sets, typically described by
probability distributions. Therefore,
 the proper mathematical context for studying models of this kind
is the theory of random point fields on $\mathbb{R}^d$, cf. \cite[page 1311]{Neu}. In this approach, the phase space
is the set of point configurations
\begin{equation}
  \label{1}
  \Gamma = \{ \gamma \subset \mathbb{R}^d: |\gamma \cap \Lambda| < \infty \  {\rm for}  \  {\rm any}  \ {\rm compact}  \ \Lambda \subset \mathbb{R}^d \},
\end{equation}
where $|\cdot|$ denotes cardinality. A proper subset of $\Gamma$ is
the set of finite configurations $\Gamma_0 = \{\gamma \in \Gamma:
|\gamma| \in \mathbb{N}_0\}$. The set $\Gamma$ can be equipped with
the {\it vague topology}, see \cite{Albev,Tobi}, and thus with the
corresponding Borel $\sigma$-field $\mathcal{B}(\Gamma)$. The
elements of $\Gamma$ are considered as point states in the sense
that, for  a suitable function $F:\Gamma\to \mathbb{R}$, the number
$F(\gamma)$ is the value of {\em observable} $F$ in state $\gamma$.
Along with such point states one employs those
determined by probability measures on $\mathcal{B}(\Gamma)$. Then
 the corresponding value is the integral
\begin{equation}
  \label{1a}
 \langle \! \langle F, \mu \rangle \! \rangle := \int_{\Gamma} F d \mu,
\end{equation}
and the system's dynamics are described as maps $\mu_0 \mapsto
\mu_t$, $t>0$.  The elementary acts here include birth, death,
immigration, jumps, diffusion, etc. In the Markov approach, the map
$\mu_0 \mapsto \mu_t$ is obtained from the Fokker-Planck equation
\begin{equation}
  \label{2}
  \frac{d}{dt} \mu_t = L^\mu \mu_t, \qquad \mu_t|_{t=0} = \mu_0, \quad t\geq 0,
\end{equation}
in which `operator' $L^\mu$ specifies the model.
By the duality
\begin{equation}
 \label{R}
\langle \! \langle F_0, \mu_t \rangle \! \rangle = \langle \!
\langle F_t, \mu_0 \rangle \! \rangle, \qquad t>0,
\end{equation}
the observed evolution $\langle \! \langle F, \mu_0
\rangle \!  \rangle \mapsto \langle \! \langle F, \mu_t \rangle \!
\rangle$
can also be considered as the evolution $\langle \! \langle F_0, \mu \rangle \!  \rangle
\mapsto \langle \! \langle F_t, \mu \rangle \!  \rangle$, obtained from the Kolmogorov equation
\begin{equation*}
\frac{d}{d t} F_t = L F_t , \qquad F_t|_{t=0} = F_0,
\end{equation*}
where $L$ and $L^\mu$ are dual in the sense of (\ref{R}). A number of such models
 are discussed in \cite{Dima2}. In this article,
we deal with the {\em spatial logistic model} specified by
\begin{eqnarray}\label{R20}
(LF)(\gamma) &=& \sum_{x\in \gamma}\left[ m + E^{-} (x, \gamma\setminus x) \right]\left[F(\gamma\setminus x)
- F(\gamma) \right]\\[.2cm]
&& \qquad + \int_{\mathbb{R}^d} E^{+} (y, \gamma) \left[ F(\gamma\cup y) - F(\gamma) \right]dy, \nonumber
\end{eqnarray}
where
\begin{equation}
 \label{Ra20}
E^{+}(x, \gamma) = \sum_{y\in \gamma} a^{+} (x- y).
\end{equation}
The first term in (\ref{R20}) describes the death of the particle
located at $x$ occurring independently with rate $m\geq 0$
(intrinsic mortality), and under the influence of the other
particles in $\gamma$ with rate $E^{-}(x, \gamma\setminus x)\geq 0$
(competition). The second term in (\ref{R20}) describes the birth of
a particle at $y\in \mathbb{R}^d$ given by the whole configuration
$\gamma$ with rate $E^{+}(y, \gamma)\geq 0$. Models of this kind
appear in population biology, see \cite{FM,BCFKKO} for more detail.
 A particular case of (\ref{R20}) with $E^{-}\equiv 0$ is  the {\em
contact model} \cite{Dima,KKP,KS}.

The sums in (\ref{R20}) and  (\ref{Ra20}) are taken over possibly
infinite sets, and thus the very definition of $L^\mu$ in (\ref{2})
gets highly problematic if one wants to include in its domain
nontrivial states with `infinite number of entities'. For the same
reason, direct solving (\ref{2}) for such states is also beyond the
possibilities of the existing mathematical methods.  For infinite
systems of physical particles, N. N. Bogoliubov suggested studying
the dynamics of states indirectly -- by means of correlation
functions,  see \cite{Dob}. Their evolution  is obtained from an
infinite chain of equations, so called {\it BBGKY hierarchy}, that
links to each other correlation functions of different order.
Starting from the late 1990'th, a similar approach is being
implemented in studying Markov dynamics on the phase space $\Gamma$,
see \cite{FKK,BCFKKO} and the references quoted therein. In this
approach, various aspects of the evolution
of the model specified by a particular case of (\ref{R20}) with `infinite number of
entities' were studied in \cite{FKKK,FKK-MFAT,FKK}. The results of
\cite{FKKK} were announced in \cite{EPJ}. The aim of this work is to
describe the evolution of states of the model (\ref{R20}) with
`finite number of entities'.

\section{The Setup and Results}
\label{sec:2}

\subsection{Preliminaries}
\label{ssec:2.1}

By $\mathcal{B}(\mathbb{R}^d)$ and $\mathcal{B}_{\rm
b}(\mathbb{R}^d)$ we denote the set of all Borel and all bounded
Borel subsets of $\mathbb{R}^d$, respectively. The set of
configurations $\Gamma$ defined in (\ref{1}) is equipped with the
vague topology -- the weakest topology that makes the
maps
\[
\Gamma \ni \gamma \mapsto \sum_{x \in \gamma} f(x) \in \mathbb{R}
\]
continuous for all compactly supported continuous functions
$f:\mathbb{R}\to \mathbb{R}$. This topology can be completely and
separably metrized that turns $\Gamma$ into a Polish spaces, see
\cite{Albev,Tobi}. By $\mathcal{B}(\Gamma)$ and
$\mathcal{P}(\Gamma)$ we denote the Borel $\sigma$-field of subsets
of $\Gamma$ and the set of all probability measures on
$\mathcal{B}(\Gamma)$, respectively.

The set of finite configurations
\begin{equation*}
  \Gamma_{0} = \bigcup_{n\in \mathbb{N}_0} \Gamma^{(n)} ,
\end{equation*}
is the disjoint union of the sets of $n$-particle configurations:
\begin{equation*}
\Gamma^{(0)} = \{ \emptyset\}, \qquad \Gamma^{(n)} = \{\gamma \in \Gamma: |\gamma| = n \}, \ \ n\in \mathbb{N}.
\end{equation*}
For $n\geq 2$, $\Gamma^{(n)}$ can be identified with the symmetrization of the set
\begin{equation}
  \label{7a}
\left\{(x_1, \dots , x_n)\in \bigl(\mathbb{R}^{d}\bigr)^n: x_i \neq x_j, \ {\rm for} \ i\neq j\right\},
\end{equation}
which allows one to introduce the corresponding (Euclidean) topology
on $\Gamma^{(n)}$. Then  the topology
on the whole $\Gamma_0$ is defined as follows:  $A\subset \Gamma_0$ is said to be
open if its intersection with each $\Gamma^{(n)}$ is open. This (intrinsic)
topology differs from that induced on $\Gamma_0$ by the vague
topology of $\Gamma$. At the same time,  $\Gamma_0\in\mathcal{B}(\Gamma)$.
Thus, a function $G :\Gamma_0 \to \mathbb{R}$ is
measurable as a function on $\Gamma$ if and only if its restrictions
to each $\Gamma^{(n)}$ are Borel functions. Clearly, these
restrictions fully determine $G$. In view of (\ref{7a}), the
restriction of $G$ to $\Gamma^{(n)}$ can be extended to a symmetric
Borel function $G^{(n)}: (\mathbb{R}^d)^n \to \mathbb{R}$, $n\in
\mathbb{N}$, such that
\begin{equation}
  \label{8}
 G(\gamma) = G^{(n)} (x_1, \dots , x_n), \qquad {\rm for} \ \  \gamma = \{x_1 , \dots , x_n\}.
\end{equation}
It is convenient to complement (\ref{8})  by putting
$G(\emptyset)=G^{(0)} \in \mathbb{R}$. Let all $G^{(n)}$, $n\in \mathbb{N}$, be
continuous functions with compact support and $G$ be as in (\ref{8}).
For all such $G$, we then write
\begin{eqnarray}
  \label{9}
\int_{\Gamma_0}  G(\gamma) \lambda (d \gamma) = G^{(0)} + \sum_{n=1}^\infty \frac{1}{n!} \int_{(\mathbb{R}^d)^n} G^{(n)} (x_1 , \dots x_n) d x_1 \cdots d x_n .
\end{eqnarray}
This expression determines a
$\sigma$-finite measure $\lambda$ on $\Gamma_0$, called the {\em
Lebesgue-Poisson} measure.

Having in mind that $\Gamma_0 \in \mathcal{B}(\Gamma)$ we define
\[
\mathcal{F}(\Gamma_0) = \{ B \subset \Gamma_0 : B \in \mathcal{B}(\Gamma)\}.
\]
The intrinsic topology of $\Gamma_0$ is clearly metrizable. The corresponding metric space is also complete and separable. The embedding $\Gamma_0 \hookrightarrow \Gamma$ is continuous and hence $\mathcal{B}(\Gamma_0)/\mathcal{B}(\Gamma)$-measurable. By the Kuratowski theorem \cite[page 21]{Pa} we then have
\begin{equation*}
\mathcal{F}(\Gamma_0) = \mathcal{B}(\Gamma_0).
\end{equation*}
Thus, each $\mu\in \mathcal{P}(\Gamma)$ with the property
\begin{equation}
  \label{Kur1}
 \mu(\Gamma_0) =1
\end{equation}
can be redefined as a probability measure on
$\mathcal{B}(\Gamma_0)$.

The considered system of entities is finite if it is in a state
  with property (\ref{Kur1}).
In the probabilistic interpretation, such a system has a `random
number of entities', whereas $\mu(\Gamma^{(n)})$ is the probability
that this number is $n$. Correspondingly,
\begin{equation}
  \label{Kur2}
N_\mu:=\int_{\Gamma_0} |\gamma| \mu( d \gamma) = \sum_{n=0}^\infty n
\mu(\Gamma^{(n)})
\end{equation}
is the expected number of entities in state $\mu$. Of course, the
convergence in (\ref{Kur2}) depends on the asymptotic properties of
the sequence $\{\mu(\Gamma^{(n)})\}_{n\in \mathbb{N}_0}$, which
characterize the state. For the Poisson measure $\pi_\varrho$ with
density $\varrho$, we have
\begin{equation*}
N_{\pi_\varrho} = \int_{\mathbb{R}^d} \varrho (x) d x,
\end{equation*}
and
\begin{equation*}
 \pi_\varrho(\Gamma^{(n)}) =  N_{\pi_\varrho}^n \exp \left( - N_{\pi_\varrho} \right)/n!,
\end{equation*}
which yields
\begin{equation}
  \label{Kur4}
 \int_{\Gamma_0} \exp\left( \beta |\gamma| \right) \pi_\varrho ( d \gamma) = \exp\left( (e^\beta -1)\int_{\mathbb{R}^d} \varrho (x) d x \right) ,
\end{equation}
holding for all $\beta \geq 0$. Another relevant property of
$\pi_\varrho$ is that it is absolutely continuous with respect to
the Lebesgue-Poisson measure introduced in (\ref{9}). Furthermore,
\begin{equation}
  \label{Ku1}
  \pi_\varrho ( d \gamma) = \left( \prod_{x\in \gamma}
  \varrho(x)\right) \exp \left( - N_{\pi_\varrho} \right) \lambda (
  d \gamma).
\end{equation}
The aim of this work is to show that the evolution described by
(\ref{2}) and (\ref{R20}) can preserve the properties of having finite
exponential moments as in (\ref{Kur4}) and of being absolutely
continuous with respect to the Lebesgue-Poisson measure. We do this in Theorem \ref{1tm} below.

\subsection{The results}
\label{eveqs}

For finite systems, one can restrict the sums in (\ref{R20}) and
(\ref{Ra20}) to finite $\gamma$, which in turn allows for writing
the Fokker-Planck equation (\ref{2}) in the form
\begin{eqnarray}
  \label{FPe}
\frac{d}{dt} \mu_t (d \eta) & = & (L^\mu \mu)(d \eta) \\[.2cm]& = & - \Xi (\mathbb{R}^d, \eta) \mu_t (d
\eta) + \int_{\Gamma_0} \Xi ( d \eta, \gamma) \mu_t ( d \gamma),
\quad \mu_t|_{t=0} = \mu_0, \nonumber
\end{eqnarray}
where, for $A \in \mathcal{B}(\Gamma_0)$ and $\gamma\in \Gamma_0$,
\begin{equation}
  \label{FP1}
 \Xi (A, \gamma) = \sum_{x\in \gamma} ( m + E^{-} (x , \gamma \setminus x)) \mathbb{I}_A (\gamma \setminus x) + \int_{\mathbb{R}^d}
 E^{+} (y , \gamma )) \mathbb{I}_A (\gamma \cup y) dy,
\end{equation}
 which is a
measure kernel on $(\Gamma_0, \mathcal{B}(\Gamma_0))$. The kernels
in (\ref{Ra20}), and hence in (\ref{FP1}), are supposed to satisfy:
(i) $E^{-}$ is nonnegative and measurable in each of its arguments;
 (ii) $a^{+}(x) = a^{+}(-x) \geq 0$ and
\begin{equation}
 \label{AA}
\qquad a^{+} \in L^1(\mathbb{R}^d) .
\end{equation}
Set
\begin{equation}
  \label{AA1}
\langle a^{+} \rangle = \int_{\mathbb{R}^d} a^{+} (x) d x, \qquad   E^{-} (\gamma) = \sum_{x\in \gamma} E^{-}(x, \gamma \setminus x).
\end{equation}
Then
\begin{equation}
  \label{AA2}
0\leq \Xi (\mathbb{R}^d, \gamma)= m|\gamma| + E^{-} (\gamma) + \langle a^{+} \rangle |\gamma| < \infty, \qquad \gamma \in \Gamma_0.
\end{equation}
Now let $\mathcal{M}$ be the set of all finite signed measures
(i.e., $\sigma$-additive functions) on $\mathcal{B}(\Gamma_0)$, and
let $\mathcal{M}_{+}$ be its subset consisting of nonnegative
measures. By means of the Jordan decomposition $\mu = \mu_{+} -
\mu_{-}$, $\mu_{\pm} \in \mathcal{M}_{+}$, we introduce the norm
\begin{equation}
  \label{FP2}
 \|\mu\| = \mu_{+} (\Gamma_0) + \mu_{-} (\Gamma_0),
\end{equation}
which makes $\mathcal{M}$ a Banach space. The right-hand side of
(\ref{FPe})  defines an unbounded linear operator $L^\mu:
\mathcal{M}\to \mathcal{M}$ with domain
\begin{equation}
  \label{FP3}
  \mathcal{D} (L^\mu) = \{ \mu \in \mathcal{M}: \Xi(\mathbb{R}^d, \cdot) \mu_{\pm} \in \mathcal{M}\}.
\end{equation}
\begin{definition}
\label{1zdf} By a classical solution of the problem in
(\ref{FPe}), on the time interval $[0, \tau)$, $\tau \leq +\infty$ and in the space $\mathcal{M}$,  we understand a map $[0,\tau)\ni t
\mapsto \mu_{t} \in \mathcal{D} (L^\mu) $ which is continuously
differentiable on $[0,\tau)$ and such that (\ref{FPe}) is
satisfied  for all $t\geq 0$. Such a solution is said to be global if $\tau= +\infty$.
\end{definition}
Let us turn now to solving (\ref{FPe}) in the class of measures
absolutely continuous with respect to the Lebesgue-Poisson measure
$\lambda$. To this end we set
\begin{equation}
  \label{FP30}
  \mu_t ( d \gamma) = R_t (\gamma) \lambda (d \gamma).
\end{equation}
By means of the following standard relation
\begin{equation*}
 \int_{\Gamma_0} \left( \sum_{x\in \gamma} f(x, \gamma\setminus x)\right)
  \lambda (d \gamma) = \int_{\Gamma_0} \int_{\mathbb{R}^d}  f(x, \gamma)
  \lambda (d \gamma) d x,
\end{equation*}
see e.g. \cite[Lemma 2.1]{FKK}, one transforms (\ref{FPe}) into the
problem
\begin{eqnarray}
  \label{FPd}
 \frac{d}{dt} R_t (\gamma) & = & -\left[\sum_{x\in \gamma}(m+ \langle a^{+} \rangle + E^{-}(x,\gamma \setminus x))\right] R_t (\gamma) \\[.2cm]
&  + & \int_{\mathbb{R}^d} (m+E^{-}(x,\gamma)) R_t (\gamma \cup x) d
x + \sum_{x\in \gamma} E^{+} (x, \gamma \setminus x) R_t (\gamma
\setminus x), \nonumber \\[.2cm]
& & R_t|_{t=0}  =  R_0. \nonumber
\end{eqnarray}
Correspondingly, the solution of (\ref{FPd}) is to be sought in the
space $L^1(\Gamma_0, \lambda)$.

Let $\mathcal{M}_{+,1}$ denote the set of $\mu\in \mathcal{M}$ such
that $\mu\geq 0$ and $\|\mu\|=1$. That is, $\mathcal{M}_{+,1}$ is
the set of all probability measures on $\mathcal{B}(\Gamma_0)$. We
shall also use the following sets
\begin{eqnarray}
  \label{Q1}
 \mathcal{M}^{(n)} & = & \bigg{\{} \mu \in \mathcal{M}: \int_{\Gamma_0}
|\eta|^n \mu_{\pm }(d \eta) < \infty\bigg{\}}, \qquad \ \quad n\in
\mathbb{N}, \\[.2cm]
 \mathcal{M}^\beta & = & \bigg{\{} \mu \in \mathcal{M}: \int_{\Gamma_0}
\exp ( \beta|\eta|) \mu_{\pm }(d \eta) < \infty\bigg{\}}, \quad
\beta >0. \nonumber
\end{eqnarray}
Finally, we denote
\[
\mathcal{M}^{(n)}_{+,1} = \mathcal{M}^{(n)} \cap \mathcal{M}_{+,1},
\qquad \mathcal{M}^{\beta}_{+,1} = \mathcal{M}^{\beta} \cap
\mathcal{M}_{+,1},
\]
and also
\begin{equation}
  \label{Q2}
 T(\alpha, \beta) = (\alpha - \beta) e^{-\alpha}/ \langle a ^{+}
 \rangle, \qquad 0 < \beta < \alpha.
\end{equation}
\begin{theorem}
  \label{1tm}
The problem in (\ref{FPe}) with $\mu_0 \in \mathcal{M}^{(1)}$ has a
unique global classical solution $\mu_t\in \mathcal{M}$ which has
the following properties: 
\begin{itemize}
\item[(a)] \ \ for each $n\in \mathbb{N}$,  $\mu_0 \in
\mathcal{M}^{(n)}_{+,1}$ implies $\mu_t \in \mathcal{M}^{(n)}_{+,1}$
for all $t>0$; \vskip.1cm
\item[(b)] \ \ if the parameters in (\ref{FP1}) satisfy:
(1) $\langle a^{+} \rangle =0$, or (2) $m > \langle a^{+} \rangle >0$,
then $\mu_0 \in \mathcal{M}^{\beta}_{+,1}$ implies $\mu_t \in
\mathcal{M}^{\beta}_{+,1}$  for all $t>0$, which holds for all
positive $\beta$ in case (1) and all $\beta \in (0, \log m -
\log \langle a^{+} \rangle)$ in case (2); \vskip.1cm
\item[(c)] \ \ for each $\beta^*>0$ and $\beta_* \in (0,
\beta^*)$, $\mu_0 \in \mathcal{M}^{\beta^*}_{+,1}$ implies $\mu_t
\in \mathcal{M}^{\beta_*}_{+,1}$  for all $t\in (0,T_*)$, where $T_*
= T(\beta^*, \beta_*)$, see (\ref{Q2}); \vskip.1cm
\item[(d)] \ \ for all $t>0$, $\mu_t$ has the property (\ref{FP30})
if $\mu_0$ does so.
\end{itemize}
\end{theorem}
Let us make some comments on the above statement. Under more
restrictive assumptions, part {\it (a)} was obtained in \cite{FM}, see
Assumption A and Theorem 3.1, respectively. The properties claimed
in {\it (b) -- (d)} resemble those for the Poisson measure, see
(\ref{Kur4}) and (\ref{Ku1}). Note, however, that only in case {\it (c1)} we
have the existence of the exponential moments for {\it all}
$\beta>0$. Also note that the only property of $E^{-}$ in (\ref{FP1}) which we need is its measurability, and that
the above results hold for $E^{-}\equiv 0$. That is, the competition plays no role if the system is finite, unlike for
the same model with infinite number of entities studied in \cite{EPJ,FKKK}.

\section{Proof of Theorem \ref{1tm}}

The proof is crucially based on the perturbation theory for
semigroups of positive operators developed in \cite{TV}. In the next
subsection, we present its relevant aspects.

\subsection{The Thieme-Voigt perturbation theory}

Let $\Psi$ be either $\mathcal{M}$ or $L^1(\Gamma_0, \lambda)$. In
both cases, the norm $\|\cdot\|_\Psi$ is additive on the cone of
positive elements $\Psi_{+}$. By $\varphi$ we denote the positive
linear functional on $\Psi$ such that $\varphi (\psi) =
\|\psi\|_\Psi$ whenever $\psi \in \Psi_{+}$. A semigroup $S=
\{S(t)\}_{t\geq 0}$ of linear operators on $\Psi$ is called
substochastic (respectively, stochastic) if: (a) it is strongly
continuous in $t$ ($C_0$-property); (b) $S(t):\Psi_{+} \to \Psi_{+}$
for all $t\geq 0$ (positivity); (c) $\|S(t)\psi \|_\Psi \leq
\|\psi\|_{\Psi}$ (respectively, $\|S(t)\psi \|_\Psi =
\|\psi\|_{\Psi}$) for all $t\geq 0$ and $\psi \in
\Psi_{+}$. The next statement is the
relevant part of Theorem 2.2 of \cite{TV}.
\begin{proposition}
  \label{1pn}
Let $S_0$ be a positive $C_0$-semigroup on $\Psi$, with generator
$(A_0, \mathcal{D}(A_0))$. Let also $B: \mathcal{D}(A_0) \to \Psi$
be a positive linear operator such that
\begin{equation}
  \label{B}
\varphi((A_0+ B)\psi) \leq 0, \qquad \psi \in \mathcal{D}(A_0) \cap
\Psi_{+}.
\end{equation}
Then for all $r\in (0,1)$, the operator $(A_0 + rB,
\mathcal{D}(A_0))$ is the generator of a substochastic semigroup on
$\Psi$.
\end{proposition}
Assume that there exists a subspace  $\widetilde{\Psi} \subset \Psi$
such that: (a) there exists a norm $\|\cdot
\|_{\widetilde{\Psi}}$ and a positive linear functional
$\tilde{\varphi}$ on $\widetilde{\Psi}$ such that $\|\psi\|_{\widetilde{\Psi}} =
\tilde{\varphi}(\psi)$ for $\psi\in \widetilde{\Psi}_{+} :=
\widetilde{\Psi} \cap \Psi_{+}$; (b) $(\widetilde{\Psi}, \|\cdot
\|_{\widetilde{\Psi}})$ is a Banach space; (c)
the embedding $\widetilde{\Psi} \hookrightarrow \Psi$
is continuous. Let $S_0$ be as in
Proposition \ref{1pn}. Assume that each $S_0(t)$, $t>0$, leaves
$\widetilde{\Psi}$ invariant and the restrictions
$\widetilde{S}_0(t)$ of $S_0(t)$ to ${\widetilde{\Psi}}$ constitute
a $C_0$-semigroup $\widetilde{S}_0$ on $\widetilde{\Psi}$, generated
by the restriction of $A_0$ to
\begin{equation}
  \label{B1}
\mathcal{D}(\widetilde{A}_0) = \{ \psi\in \mathcal{D}(A_0) \cap
\widetilde{\Psi}: A_0 \psi \in \widetilde{\Psi}\}.
\end{equation}
For $B$ as in Proposition \ref{1pn}, assume additionally that
$B:\mathcal{D}(\widetilde{A}_0) \to \widetilde{\Psi}$. The next
statement is the relevant part of Theorem 2.7 of \cite{TV}.
\begin{proposition}
\label{2pn} Let the assumptions above including (\ref{B}) be
satisfied. Additionally, assume that $-A_0$ is positive and there
exists positive $c$ and $\varepsilon$ such that the following holds
\begin{equation}
  \label{B2}
\tilde{\varphi} ((A_0+B)\psi) \leq c\|\psi\|_{\widetilde{\Psi}} -
\varepsilon \|A_0 \psi\|_\Psi, \qquad \psi \in
\mathcal{D}(\widetilde{A}_0) \cap \Psi_{+}.
\end{equation}
Then the closure  of $(A_0+B, \mathcal{D}(A_0))$ in $\Psi$ is the
generator of a substochastic semigroup $S$ in $\Psi$ that leaves
$\widetilde{\Psi}$ invariant. If the equality in (\ref{B}) holds,
then $S$ is stochastic.
\end{proposition}

\subsection{The proof}

\underline{Claim (a)}
\vskip.1cm \noindent
Set $\Psi = \mathcal{M}$, $\widetilde{\Psi}= \mathcal{M}^{(n)}$, see (\ref{Q1}), and
\begin{equation}
  \label{B3}
\widetilde{\varphi} (\mu) = \int_{\Gamma_0} |\gamma|^n \mu( d
\gamma).
\end{equation}
Note that $\varphi(\mu) = \mu(\Gamma_0)$, cf. (\ref{FP2}). Next, set
\begin{equation}
  \label{B4}
(A_0 \mu) (d \gamma) = - \Xi(\mathbb{R}^d, \gamma) \mu( d \gamma),
\qquad \mathcal{D}(A_0) = \mathcal{D}(L^\mu),
\end{equation}
where the latter set is defined in (\ref{FP3}). The operator $(A_0, \mathcal{D}(A_0))$
generates the semigroup $S_0$ on $\mathcal{M}$ defined by
\begin{equation*}
(S_0 (t) \mu (d \gamma) = \exp \left( - t \Xi(\mathbb{R}^d, \gamma)
\right)  \mu( d \gamma),
\end{equation*}
which is clearly substochastic. Then the operator, cf. (\ref{FPe}),
\begin{equation}
  \label{B6}
(B\mu)(d \gamma) = \int_{\Gamma_0} \Xi(d \gamma, \eta) \mu( d
\eta)
\end{equation}
is positive and the following holds, cf. (\ref{FP2}) and (\ref{B4}),
\[
\|B\mu\| = \|A_0\mu\| = \int_{\Gamma_0}\Xi(\mathbb{R}^d, \gamma)
\mu_{+}( d \gamma) + \int_{\Gamma_0}\Xi(\mathbb{R}^d, \gamma)
\mu_{-}( d \gamma),
\]
that is, $B: \mathcal{D}(A_0) \to \mathcal{M}$. Moreover, for
positive $\mu \in \mathcal{D}(A_0)$, we have
\begin{eqnarray*}
\varphi((A_0 + B)\mu) & = & - \int_{\Gamma_0}\Xi(\mathbb{R}^d,
\gamma) \mu( d \gamma) + \int_{\Gamma_0}\Xi(\mathbb{R}^d, \gamma)
\mu( d
\gamma) = 0.
\end{eqnarray*}
The set defined in (\ref{B1}) consists of those $\mu\in
\mathcal{M}^{(n)}$ for which
\begin{equation}
  \label{B7a}
\int_{\Gamma_0} |\gamma|^n \Xi(\mathbb{R}^d, \gamma) \mu_{\pm} ( d
\gamma) < \infty.
\end{equation}
Then $B$ defined in (\ref{B6}) maps this set into
$\mathcal{M}^{(n)}$. Hence, to use Proposition \ref{2pn} it is left to check the validity of the
corresponding version of (\ref{B2}), the left-hand side
of which can be written in the form, cf. (\ref{B3}) and (\ref{1a}),
\begin{eqnarray*}
\tilde{\varphi}((A_0 + B)\mu)  & = & \int_{\Gamma_0} F (\gamma) (L^\mu \mu)(d \gamma) =
\langle \! \langle F L^\mu \mu \rangle \! \rangle \\[.2cm] & = & \langle \! \langle LF \mu \rangle \! \rangle
= \int_{\Gamma_0}(L F) (\gamma) \mu(d \gamma),   \nonumber
\end{eqnarray*}
where $F(\gamma) = |\gamma|^n$ and $L$ is given in (\ref{R20}). Thus,
(\ref{B2}) now takes the form
\begin{equation*}
\int_{\Gamma_0} \bigg{( }(LF)(\gamma) - c F(\gamma) + \varepsilon \Xi (\mathbb{R}^d, \gamma) \bigg{)} \mu( d \gamma) \leq 0
\end{equation*}
which has to hold for all positive $\mu$ with the property (\ref{B7a}), including $\mu= \delta_\eta$ for each $\eta \in \Gamma_0$.
This amounts to the following
\begin{eqnarray}
  \label{B10}
& & - \sum_{x\in \gamma} \left[ m + E^{-} (x, \gamma\setminus x)\right]\left(F(\gamma) - F(\gamma\setminus x) - \varepsilon \right) \\[.2cm]
& & \quad + \int_{\mathbb{R}^d}  E^{+} (y, \gamma)\left(F(\gamma\cup y) - F(\gamma) + \varepsilon \right)d y \leq c F(\gamma),\nonumber
\end{eqnarray}
which has to hold for some positive $c$ and $\varepsilon$ and all $\gamma \in \Gamma_0$. For $F(\gamma) = |\gamma|^n$ and $\varepsilon =1$,
we have that $F(\gamma) - F(\gamma\setminus x) - \varepsilon \geq 0$ for all $\gamma\neq \emptyset$. Then (\ref{B10}) is satisfied if
the following holds, cf. (\ref{AA}),
\[
\langle a^{+} \rangle |\gamma| \left[(|\gamma| + 1)^n - |\gamma|^n + 1 \right] \leq c |\gamma|^n,
\]
which is the case for $c= \langle a^{+} \rangle 2^{n+1}$. Thus, by Proposition \ref{2pn} the closure of $(A_0+ B, \mathcal{D}(A_0))$ as given in (\ref{B4}), (\ref{B6}), and
(\ref{FP3}) generates a stochastic semigroup $S$ on $\mathcal{M}$, which leaves invariant $\mathcal{M}^{(n)}$.
Then the solution in question is obtained as $\mu_t = S(t) \mu_0$.
\vskip.1cm \noindent
\underline{Claim (b)}
\vskip.1cm \noindent
In this case, we set $\widetilde{\Psi} = \mathcal{M}^\beta$ and find $\beta>0$ for which (\ref{B10})
holds with $F(\gamma) = \exp( \beta |\gamma|)$. For $\langle a^{+} \rangle = 0$, we fix arbitrary $\beta >0$ and then take $\varepsilon \leq
e^\beta -1$. In this  case, (\ref{B10}) holds for each $c>0$ and $\gamma \in \Gamma_0$. For $\langle a^{+} \rangle >0$, (\ref{B10}) holds, with each $c>0$, for positive $\beta$ and $\varepsilon$ that satisfy
\begin{equation}
\label{B11}
\langle a^{+} \rangle e^{2\beta} - (\langle a^{+} \rangle + m) e^\beta + m =: P(e^\beta) \leq - \varepsilon (\langle a^{+} \rangle + m).
\end{equation}
The polynomial $P$ has two roots $1$ and $m / \langle a^{+} \rangle$. Thus, for each fixed $\beta$ such that the following holds $1< e^\beta < m / \langle a^{+} \rangle$, one finds $\varepsilon$ for which
also (\ref{B11}) holds. This completes the proof  in this case.

\vskip.1cm \noindent
\underline{Claim (c)}
\vskip.1cm \noindent
Fix $\beta^*>0$ such that $\mu_0 \in \mathcal{M}_{+,1}^{\beta^*}$. By claim (a) we already know that $\mu_t = S(t) \mu_0$ is
the unique global solution of (\ref{FPe}) in $\mathcal{M}$. For a given $\beta_*\in (0,\beta^*)$, let us show that (\ref{FPe}) has also a classical solution in
$\mathcal{M}^{\beta_*}$ on the time interval $(0,T_*)$. Since
$\mathcal{M}^{\beta_*}\subset \mathcal{M}$, this solution would coincide with $S(t) \mu_0$ and hence
be in $\mathcal{M}_{+,1}^{\beta_*}$.

Thus, we consider (\ref{FPe}) in $\mathcal{M}^{\beta_*}$ with $\mu_0 \in \mathcal{M}_{+,1}^{\beta^*}$.
We study this problem by means of the scale $\{\mathcal{M}^\beta: \beta \in [\beta_*, \beta^*]\}$ of Banach spaces.
Note that the embedding $\mathcal{M}^\beta \hookrightarrow \mathcal{M}^{\beta'}$, $\beta > \beta'$, is dense and continuous.
For such $\beta$, we define
$\mathcal{M}^{\beta}$ by
\begin{equation}
  \label{C1}
 \varphi_\beta(\mu) = \int_{\Gamma_0} \exp(\beta |\gamma|) \mu( d \gamma).
\end{equation}
Let also $A_0$ be defined on $\mathcal{M}^{\beta}$ by (\ref{B4}) with domain
\begin{equation}
  \label{C2}
\mathcal{D}_{\beta} (A_0) = \{ \mu \in \mathcal{M}^{\beta}: \Xi (\mathbb{R}^d, \cdot) \mu \in \mathcal{M}^{\beta}\}.
\end{equation}
We split the remaining part of $L^\mu$ in (\ref{FPe}) into two terms, $B$ and $C$, defined on $\mathcal{M}^{\beta}$ as follows
\begin{equation}
  \label{C3}
(B\mu)(A) = \int_{\Gamma_0}\left( \sum_{x\in \gamma} \left[ m + E^{-} (x, \gamma\setminus x)\right]\mathbb{I}_A (\gamma\setminus x)\right) \mu( d \gamma),
\end{equation}
and
\begin{equation}
  \label{C3a}
(C\mu)(A) = \int_{\Gamma_0} \left(\int_{\mathbb{R}^d} E^{+} (y ,\gamma) \mathbb{I}_A ( \gamma \cup y) dy \right) \mu(d \gamma) \quad A \in \mathcal{B}(\Gamma_0).
\end{equation}
Let us show that $B: \mathcal{D}_{\beta} (A_0) \to \mathcal{M}^\beta$. To this end we take $\mu\in \mathcal{D}_{\beta} (A_0)\cap \mathcal{M}_{+}$ and calculate, see (\ref{C3}), (\ref{AA1}), and (\ref{AA2}),
\begin{eqnarray}
\label{C4}
\|B \mu\|_{\mathcal{M}^\beta}& = & \varphi_\beta ( B\mu)\\[.2cm]& = & \int_{\Gamma_0} \exp(\beta |\eta|)
\int_{\Gamma_0} \left( \sum_{x\in \gamma} \left[ m + E^{-} (x, \gamma\setminus x)\right]\delta_{\eta \setminus x} ( d \eta) \right)\mu (d\gamma)  \nonumber \\[.2cm]
& = & e^{-\beta} \int_{\Gamma_0}\exp(\beta |\gamma|) \left(m |\gamma| + E^{-}(\gamma) \right) \mu( d \gamma) \nonumber \\[.2cm]
& \leq & e^{-\beta} \int_{\Gamma_0} \exp(\beta |\gamma|) \Xi(\mathbb{R}^d, \gamma)\mu( d \gamma) \nonumber.
\end{eqnarray}
The latter estimate and (\ref{B4}) and (\ref{C1}) yield
\begin{equation*}
 \varphi_\beta ((A_0 + r B)\mu) \leq 0, \qquad \mu \in \mathcal{D}_\beta (A_0) \cap \mathcal{M}_{+},
\end{equation*}
holding for some $r\in (0,1)$.
Hence, by Proposition \ref{1pn} $((A_0 + B), \mathcal{D}_\beta (A_0))$ is the generator of a substochastic semigroup $S_\beta$ on $\mathcal{M}^\beta$.
 Now for  $\beta' \in (\beta_*, \beta$, let us show that $C$ as given in (\ref{C3a}) acts as a bounded linear operator $C: \mathcal{M}^\beta  \to \mathcal{M}^{\beta'}$. To this end, as in (\ref{C4}) we take
$\mu\in \mathcal{D}_{\beta} (A_0)\cap \mathcal{M}_{+}$ and calculate, see (\ref{AA1}),
\begin{eqnarray}
  \label{C5a}
\|C \mu\|_{\mathcal{M}^{\beta'}}& = & \varphi_{\beta'} ( C \mu)\\[.2cm]& = & \int_{\Gamma_0}  \exp(\beta' |\eta|)
\int_{\Gamma_0} \left( \int_{\mathbb{R}^d} E^{+} (y ,\gamma) \delta_{\gamma \cup y} ( d \eta) d y\right)\mu (d \gamma) \nonumber \\[.2cm]
& = & e^{\beta'} \langle a^{+} \rangle \int_{\Gamma_0}\exp(\beta' |\gamma|) |\gamma| \mu( d \gamma) \nonumber \\[.2cm]
& \leq & \frac{e^{\beta'} \langle a^{+} \rangle }{e(\beta - \beta')} \int_{\Gamma_0}\exp(\beta |\gamma|)  \mu( d \gamma) = \frac{e^\beta \langle a^{+} \rangle }{e(\beta - \beta')} \|\mu\|_{\mathcal{M}^\beta}.\nonumber
  \end{eqnarray}
That is, the operator norm $\|C\|_{\beta, \beta'}$ of $C: \mathcal{M}^\beta  \to \mathcal{M}^{\beta'}$ satisfies
\begin{equation}
 \label{C5b}
  \|C\|_{\beta, \beta'} \leq \frac{e^\beta \langle a^{+} \rangle }{e(\beta - \beta')}.
\end{equation}
In obtaining the last line in (\ref{C5a}), we have used the following obvious inequality $$\exp\left[-(\beta -\beta') N \right] N \leq \frac{1} {e(\beta - \beta')}, \qquad N\in \mathbb{N}.$$
For fixed $n\in \mathbb{N}$ and $\beta \in (\beta_*, \beta^*)$, set
\begin{gather}
  \label{C7}
\mu_t^{(n)}  =  S_\beta (t) \mu_0 + \sum_{l=1}^n \int_0^t\int_0^{t_1} \cdots \int_0^{t_{l-1}} T_l (t, t_1, \dots, t_l) \mu_0 d t_l \cdots d t_1 \quad \\[.2cm]
 T_l (t, t_1, \dots, t_l) = S_{\beta_l} (t-t_1) C_{\beta_l} S_{\beta_{l-1}} (t_1-t_2) \cdots S_{\beta_l}(t_{l-1} - t_l) C_{\beta_1} S_{\beta_0}(t_l),   \nonumber
\end{gather}
where, for a fixed $l\geq1$,
\begin{equation}
  \label{C6}
\beta_s = \beta^* - s \epsilon, \qquad \epsilon  = (\beta^* - \beta)/l, \quad s = 0, \dots , l,
\end{equation}
and hence $\beta_0 = \beta^*$ and $\beta_l= \beta$. For $s=0, \dots , l$, the operator
$C_{\beta_s}$  acts from $\mathcal{M}^{\beta_{s-1}}$ to $\mathcal{M}^{\beta_{s}}$ as a bounded operator, see (\ref{C5a}) and (\ref{C5b}).
Note  that $\mu_t^{(n)} \in \mathcal{M}^\beta$ for each $t>0$ and $n \in \mathbb{N}$. Note also that each $T_l$ is continuous in $t$ and each of $t_s$,  $s= 1, \dots , l$, which follows from the strong continuity of each $S_{\beta_s}$. This yields that $\mu_t^{(n)}$ is continuously differentiable in $\mathcal{M}^\beta$. Furthermore, for each $\beta' \in [\beta_*, \beta)$, we have that the derivative as an element of $\mathcal{M}^{\beta'}$ satisfies
\begin{equation}
  \label{C6a}
 \frac{d }{dt} \mu_t^{(n)} = (A_0 + B)\mu_t^{(n)} + C\mu_t^{(n-1)},
\end{equation}
where all these operators act in $\mathcal{M}^{\beta'}$ with domains containing $\mathcal{M}^\beta$, see (\ref{C2}), (\ref{C5a}), and (\ref{C5b}).

By (\ref{C7}) and (\ref{C6}) we have
\begin{gather}
\label{C8}
\|\mu_t^{(n)} - \mu_t^{(n-1)}  \|_{\mathcal{M}^\beta} \leq   \int_0^t\int_0^{t_1} \cdots \int_0^{t_{l-1}}
\|T_l (t, t_1, \dots, t_l)\mu_0\|_{\mathcal{M^\beta}}  d t_l \cdots d t_1 \nonumber \\[.2cm]
\leq \frac{t^n}{n!} \|C_{\beta_n}C_{\beta_{n-1}} \cdots C_{\beta_1}\mu_0\|_\beta \\[.2cm]
\leq \frac{t^n}{n!} \left(\prod_{s=1}^n \|C_{\beta_s}\|_{\beta_{s-1},\beta_s} \right) \|\mu_0\|_{\beta^*} \nonumber \\[.2cm]
\leq \frac{1}{n!} \left( \frac{n}{e}\right)^n \left(\frac{te^{\beta^*} \langle a^{+} \rangle}{\beta^* - \beta} \right)^n. \nonumber
\end{gather}
By (\ref{C8}) the sequence $\{\mu_t^{(n)}\}_{n\in \mathbb{N}}$ converges in $\mathcal{M}^\beta$ to a certain $\mu_t$, uniformly on each $[0,\theta]$, $\theta < T(\beta^*, \beta)$. On the other hand,
\begin{eqnarray*}
& & \sup_{t\in [0,\theta]}  \|(A_0 + B)[\mu_t^{(n)} - \mu_t^{(n-1)}]\|_{\mathcal{M}^{\beta'}}\\[.2cm]& & \quad \qquad  \leq \left(\|A_0\|_{\beta, \beta'} + \|B\|_{\beta, \beta'} \right)\sup_{t\in [0,\theta]}\|\mu_t^{(n)} - \mu_t^{(n-1)}\|_{\mathcal{M}^\beta} , \\[.2cm]
& & \sup_{t\in [0,\theta]}  \|C[\mu_t^{(n)} - \mu_t^{(n-1)}]\|_{\mathcal{M}^{\beta'}} \leq \|C\|_{\beta, \beta'} \sup_{t\in [0,\theta]}\|\mu_t^{(n)} - \mu_t^{(n-1)}\|_{\mathcal{M}^\beta},
\end{eqnarray*}
where the operator norms can be estimated as in (\ref{C5b}). Hence, by (\ref{C8}) $\{ d \mu_t^{(n)} / dt\}_{n\in \mathbb{N}}$ converges in $\mathcal{M}^{\beta'}$, uniformly on $[0,\theta]$. Therefore, the limiting $\mu_t\in \mathcal{M}^{\beta'}$ is continuously differentiable on
each $[0, \theta]\subset [0, T(\beta^*, \beta'))$, and
\[
\mu^{(n)}_t \to \mu_t\in \mathcal{M}^{\beta'}, \qquad n \to +\infty.
\]
On the other hand, the right-hand side of (\ref{C6a}) converges in $\mathcal{M}^{\beta'}$ to $L^\mu \mu_t$. Hence, the limiting $\mu_t$
is the classical solution of (\ref{FPe}) in $\mathcal{M}^{\beta'}$ on the time interval $[0, T(\beta^*, \beta'))$. Now we set $\beta'=\beta_*$ and obtain the proof for this case.

\vskip.1cm \noindent
\underline{Claim (d)}
\vskip.1cm \noindent
We consider the problem in (\ref{FPd}) and repeat the arguments used in the proof above. By the uniqueness of the solution of (\ref{FPe})
we obtain that $\mu_t$ in Claims (a) - (c) and $\mu_t := R_t \lambda$ coincide, which yields the proof.

\vskip.2cm
\noindent
{\bf Acknowledgement:} The author thanks Yuri Kondratiev for fruitful discussions on the subject of this work.

%
%
%


\end{document}